\documentclass{amsart}

\newtheorem{Thm}{Theorem}[section]

\newtheorem{Lem}[Thm]{Lemma}

\usepackage{graphicx}
\usepackage{latexsym}
\begin{document}

\title[On strict inclusions]
      {On strict inclusions \\ in  hierarchies of convex bodies.}\

\author{Vladyslav Yaskin}

\begin{abstract}

 Let $\mathcal I_k$ be the class of convex $k$-intersection bodies in $\mathbb{R}^n$ (in the sense of
Koldobsky) and $\mathcal I_k^m$ be the class of convex origin-symmetric  bodies all of whose
$m$-dimensional central sections are $k$-intersection bodies. We show that 1)  $\mathcal
I_k^m\not\subset \mathcal I_k^{m+1}$, $k+3\le m<n$, and 2)  $\mathcal I_l \not\subset \mathcal I_k
$, $1\le k<l < n-3$.

\end{abstract}

\subjclass[2000]{52A20, 52A21, 46B04.}
\address{Department of Mathematics, University of Oklahoma, Norman, Oklahoma 73019}
 \email{vyaskin@math.ou.edu}
\thanks{The author was supported in part by the European Network PHD, FP6 Marie Curie Actions, RTN, Contract MCRN-511953.
Part of the work was done when the author was visiting Universit\'{e} de Marne-la-Vall\'{e}e.}

\maketitle

\section{Introduction}

Let $K$ and $L$ be origin-symmetric star bodies in $\mathbb{R}^n$. Following Lutwak \cite{Lu} we
say that  $K$ is  the {\it intersection body of} $L$ if the radius of $K$ in every direction is
equal to the volume of the central hyperplane section of $L$ perpendicular to this direction, i.e.
for every $\xi \in S^{n-1}$,
$$\rho_K(\xi)=\mathrm{vol}_{n-1}(L\cap\xi^\bot).$$
The closure in the radial metric of the class of intersection bodies of star bodies gives the class
of {\it intersection bodies}.


A generalization of the concept of an intersection body was introduced by Koldobsky in
\cite{K-GAFA}. Let $1\le k<n$ and let $K$ and $L$ be origin-symmetric star bodies in
$\mathbb{R}^n$. We say that $K$ is a {\it $k$-intersection body of $L$} if for every
$(n-k)$-dimensional subspace $H\subset \mathbb{R}^n$
$$\mathrm{Vol}_{k} (K\cap H^\perp ) = \mathrm{Vol}_{n-k} (L\cap H).$$
The closure in the radial metric gives the class of {\it
$k$-intersection bodies}, which will be denoted by $\mathcal I_k$.
Note that $\mathcal I_1$ is the class of intersection bodies.

Koldobsky \cite{K-Canad} introduced the concept of embedding of a normed spaces in $L_p$, $p<0$,
and in \cite{K-GAFA} he proved that $k$-intersection bodies are the unit balls of spaces that embed
in $L_{-k}$.

A well-known property of $L_p$-spaces, proved in \cite{BDK}, is that, for any $0<p<q\le 2,$ the
space $L_q$ embeds isometrically in $L_p,$ so $L_p$-spaces become larger when $p$ decreases from 2.
Koldobsky \cite{K-Canad} extended this result   to negative $p$: { Every $n$-dimensional subspace
of $L_q$, $0<q\le 2$, embeds in $L_p$ for every $-n<p<0$. }

However, it is an open problem,  whether a normed space  $X=(\mathbb{R}^n,\|\cdot\|)$ being
embedded in $L_{-p}$ for some $0<p<n-3$ implies that $X$ embeds in $L_{-q}$ for all $p<q<n$. In
particular, is it true that every $k$-intersection body is also an $m$-intersection body for $1< k
< m < n-3$?  Note that in some cases the above statement is known to be true. Since the product of
positive definite distributions is also positive definite, one immediately obtains  that if $X$
embeds in $L_{-p}$, $0<p<n$, and $p$ divides $q$, $p<q<n$, then  $X$ also embeds in $L_{-q};$ see
\cite{M1}.

Related to these are the questions of showing that these classes of bodies are different for
different values of $p$ and $q$.

 It was shown by Koldobsky \cite{K-PAMS1996}, that there is an $n$-dimensional ($n\ge 3$) Banach subspace
 of  $L_{1/2}$ that does not embed in $L_1$, (also a subspace of $L_{1/4}$ that does not embed in
$L_{1/2}$).  Later Borwein and his colleagues at the Center for Computational Mathematics at Simon
Fraser University   showed (by computer methods) that there is a Banach space that embeds in
$L_{a/64}$ but not in $L_{(a+1)/64}$ for $a=1,2$,..., $63$. Another construction was given by
Kalton and Koldobsky in \cite{KK} and allowed to extend these results to all $0<p<q\le 1$.
Schlieper \cite{Schl} used the construction from \cite{K-PAMS1996} to show that there is a normed
space that embeds in $L_{-4}$ but does not embed in $L_{-2}$, and also a normed space that embeds
in $L_{-1/3}$ but not in $L_{-1/6}$.

In this paper we extend Schlieper's result to arbitrary integers, namely, we construct examples of
origin-symmetric convex bodies which are $k$-intersection bodies, but not $l$-intersection bodies
for $1\le l<k<n-3$. We should remark that all origin-symmetric convex bodies are $k$-intersection
bodies for $k=n-1$, $n-2$ and $n-3$, see \cite[p. 78]{K-book}.

Another result that we present here is motivated by papers of Weil \cite{W}, Neyman \cite{N} and
Yaskina \cite{Ya}. Weil constructed a convex body in $\mathbb{R}^n$ ($n \ge 3$) that is not a
zonoid but all its projections onto hyperplanes are zonoids. Neyman  showed that there are
$n$-dimensional normed spaces that do not embed in $L_p$, but all their $(n-1)$-dimensional
subspaces embed in $L_p$ for $p>0$. Yaskina constructed a body in $\mathbb{R}^n$ ($n\ge 5$), which
is not an intersection body, but all of its central hyperplane sections are intersection bodies.
(Note, that all central sections of an intersection body are necessarily  intersection bodies, see
Fallert, Goodey and Weil \cite{FGW}).

Here we generalize Yaskina's construction to prove the following. Let $\mathcal I_k^m$ be the class
of convex bodies all of whose $m$-dimensional central sections are $k$-intersection bodies. There
exists an origin-symmetric convex body $K\subset \mathbb{R}^n$, such that $K\in \mathcal I_k^m$ but
$K\not\in \mathcal I_k^{m+1}$, $k+3\le m<n$. One should compare this result with the fact (see e.g.
\cite{M1}), that all central sections of a $k$-intersection body are also $k$-intersection bodies
provided that the dimension of the sections is greater than $k$. Therefore,  $\mathcal
I_k^{m+1}\subset \mathcal I_k^{m}$.

Finally let us remark that another generalization of intersection bodies was introduced by Zhang
\cite{Zh-AJM}. These are called {\it generalized $k$-intersection bodies} by Koldobsky and {\it
$k$-Busemann-Petty bodies} by E.~Milman. See  \cite{K-GAFA}, \cite[Section 4.5]{K-book}, \cite{M1},
\cite{M2} for many interesting results explaining the relation between these different
generalizations and their connection to the lower dimensional Busemann-Petty problem.

\section{Section: $\mathcal I_k^m\not\subset \mathcal I_k^{m+1}$}

Let us start with the following criterion for $k$-intersection bodies.

\begin{Thm}(Koldobsky \cite{K-GAFA})
Let $K$ be an origin-symmetric star body in $\mathbb{R}^n,$ $1\le k <n.$  $K$ is a $k$-intersection
body if and only if the Fourier transform of $\|x\|_K^{-k}$ is a positive distribution.
\end{Thm}

The main result of this section is the following Theorem.

\begin{Thm}
Let  $k+3\le m<n$. There exists an origin-symmetric convex body $K$
that belongs to $\mathcal I_k^m$, but not to $\mathcal I_k^{m+1}$.
\end{Thm}
\noindent{\bf Proof.}  For a small $\epsilon>0$  define a body $K$
by
\begin{eqnarray*}\label{Def:K}
\|x\|_K^{-k}=|x|_2^{-k} -   2 \epsilon^{m-k} \|x\|_E^{-k}, \quad x \in \mathbb{R}^n \setminus\{0\},
\end{eqnarray*}
where $|x|_2$ is the Euclidean norm and $E$ is the ellipsoid  given
by
$$\|x\|_E^{}=  \left( x_1^2+\cdots +
x_m^2 +\frac{x_{m+1}^2+\cdots +x_n^2}{\epsilon^2}\right)^{1/2}.$$

 Since $\|x\|_E^{-1}\le  |x|_2^{-1}$, we have
that $\|x\|_K^{-1}$ is positive for a small $\varepsilon$, and  so the body $K$ is well defined.

The proof of the theorem follows from the following three lemmas.

\begin {Lem}\label{K-convex}
The body $K$ is convex for small enough $\varepsilon$.
\end{Lem}

\noindent{\bf Proof.} This is a standard perturbation argument, cf. \cite[p.96]{K-book}.  By
construction, the body $K$ is obtained by perturbing the Euclidean ball. Since the latter has
strictly positive curvature, it is enough to control the first and second derivatives of the
function $\epsilon^{m-k} \|x\|_E^{-k}$. One can see that their order is $O(\epsilon^{m-k-2})$,
which is small for small enough $\epsilon$. Therefore $K$ also has positive curvature.

\qed

 Recall that the Fourier
transform of $|x|_2^{-k}$, $0<k<n$,  equals (see \cite[p. 363]{GS})
$$(|x|_2^{-k})^\wedge(y)=C_{n,k} |y|_2^{-n+k} ,$$
where
$$C_{n,k}=\frac{ 2^{n-k}\pi^{n/2}{\Gamma\left(({n-k})/{2}\right)}}{\Gamma\left( {k}/{2}\right)}.$$

In order to compute the Fourier transform for the norms of ellipsoids, note that if $T$ is an
invertible linear transformation on $\mathbb{R}^n$, then
$$(|Tx|_2^{-k})^\wedge(y)= C_{n,k} |\det T|^{-1} |(T^*)^{-1}y|_2^{-n+k}.$$

\begin{Lem}\label{V-int}
For every $m$-dimensional subspace $H$ of $\mathbb{R}^n$, the body $K \cap H$ is a $k$-intersection
body.
\end{Lem}

\noindent{\bf Proof.} We have $$\|x\|_{K\cap H} ^{-k}=|x|_{B_2\cap H}^{-k} -   2 \epsilon^{m-k}
\|x\|_{E\cap H}^ {-k}.$$ Since $E$ is an ellipsoid with semiaxes $\epsilon$ and $1$, $E\cap H$ is
also an ellipsoid with semiaxes $a_1$, ..., $a_m$ such that  $\epsilon\le a_i \le 1$, $\forall
i=1,...,m$. There is a coordinate system in $H$ such that $$\|y\|_{K\cap H}
^{-k}=\left(y_1^2+\cdots + y_{m}^2\right)^{-k/2} -   2\epsilon^{m-k}
\left(\frac{y_1^2}{a_1^2}+\cdots +\frac{y_m^2}{a_m^2}\right)^{-k/2}.$$

Taking the Fourier transform of $\|y\|_{K\cap H} ^{-k}$ in the plane $H$ we get
$$(\|y\|_{K\cap H}^{-k})^\wedge(\xi)=C_{m,k}\left(|\xi|_2^{-m+k} - 2 \epsilon^{m-k} \prod_{i=1}^m a_i   \cdot \left({a_1^2}{\xi_1^2}+\cdots
+{a_m^2}{\xi_m^2}\right)^ {(-m+k)/2}\right).$$

Let $a_j$ be the smallest semiaxis. Then for some $\lambda \ge 1$ we have $a_j=\lambda \epsilon$.
Therefore,
$$\prod_{i=1}^m a_i \le   \lambda \epsilon.$$ On the other hand if $\xi \in S^{m-1}\subset
H$, then
$$\left({a_1^2}{\xi_1^2}+\cdots +{a_m^2}{\xi_m^2}\right)^ {(-m+k)/2} \le a_j^{-m+k}= (\lambda
\epsilon)^{-m+k}.$$ Therefore,
\begin{eqnarray*} 2^{ } \epsilon^{m-k} \prod_{i=1}^m a_i   \cdot
\left({a_1^2}{\xi_1^2}+\cdots +{a_m^2}{\xi_m^2}\right)^ {(-m+k)/2}\le 2^{ } \epsilon^{m-k}
  \lambda \epsilon (\lambda \epsilon)^{-m+k} \le 2 \epsilon.
\end{eqnarray*}
So, if $\epsilon \le 1/2$, then $(\|y\|_{K\cap H}^{-k})^\wedge(\xi)\ge 0$ for all $\xi\in
S^{n-1}\cap H$ and all $H$. Therefore all $m$-dimensional sections of $K$ are $k$-intersection
bodies.

\qed

\begin{Lem}\label{K-int}
There exists an $(m+1)$-dimensional section of  $K$  which is not a $k$-intersection body.
\end{Lem}

\noindent{\bf Proof.} Let $H=\{x \in \mathbb R^n: x_{m+2}= \cdots = x_n=0\}$. Then
$$\|x\|_{K\cap H}^{-k}=\left(x_1^2+\cdots +
x_{m+1}^2\right)^{-k/2} - 2^{ } \epsilon^{m-k} \left( x_1^2+\cdots + x_m^2
+\frac{x_{m+1}^2}{\epsilon^2}\right)^{-k/2}.$$

The Fourier transform in the variables $x_1$, ..., $x_{m+1}$ equals
\begin{eqnarray*}\left(\|x\|_{K\cap H}^{-k}\right)^\wedge (\xi)& = & C _{m+1,k}\Big(\left(\xi_1^2+\cdots +
\xi_{m+1}^2\right)^{(-m+k-1)/2}- \\
&-& 2^{ } \epsilon^{m-k}   \epsilon \left( \xi_1^2+\cdots + \xi_m^2+\epsilon^2
{\xi_{m+1}^2}\right)^{(-m+k-1)/2}\Big).
\end{eqnarray*}

If $\xi = (0,...,0,1)\in S^m\subset H$, then
$$\left(\|x\|_{K\cap H}^{-k}\right)^\wedge (\xi) = C_{m+1,k}\left( 1 - 2^{ } \epsilon^{m-k}   \epsilon
\epsilon^{-m+k-1}\right) = - C_{m+1,k}  < 0.$$ Therefore $K\cap H$
is not a $k$-intersection body.

\qed

\section{Section: $\mathcal I_l \not\subset \mathcal I_k $, $l>k$}

We will need a few auxiliary lemmas.

\begin{Lem}\label{Lem:K&M}
Let $k \in \mathbb{N}\cup\{0\}$ and $f\in C^{\infty}(S^{n-1})$, $f$ is even. Let $x=(r,\theta)$ be
polar coordinates in $\mathbb{R}^n$, so that
$f(\theta)r^{-p}=f\left(\displaystyle\frac{x}{|x|_2}\right) |x|_2^{-p}$.  Then the Fourier
transform of the distribution $f(\theta)r^{-p}$, $0<p<n$, is a homogeneous degree $-n+p$ continuous
on $\mathbb{R}^n\setminus\{0\}$ function, whose values on the unit sphere can be computed as
follows.

\begin{enumerate}
\item[(i)]  If $q<2k$, $q$ is not an odd integer, then   $\forall x \in S^{n-1}$,
$$ \left(f(\theta)r^{-n+q+1}\right)^\wedge(x)=\frac{(-1)^{k+1}\pi}{2\Gamma(2k-q)\sin(\pi(2k-q-1)
/2)}$$ $$\times
\int_{S^{n-1}}|(x,\xi)|^{2k-q-1}\Delta^k\left(f(\theta)r^{-n+q+1}\right)(\xi)d\xi,$$ where $\Delta$
is the Laplace operator on $\mathbb{R}^n$.

\item[(ii)] If $q$ is an even integer, $q=2k$, then $\forall x\in S^{n-1}$,
$$\left(f(\theta)r^{-n+2k+1}\right)^\wedge(x)=(-1)^k\pi   \int_{S^{n-1}\cap x^\perp}\Delta^k\left(f(\theta)r^{-n+2k+1})\right)(\xi)d\xi. $$

\item[(iii)] If $q$  is an odd integer, $q=2k-1\ge 1 $, then $\forall x\in S^{n-1}$,
\begin{eqnarray*}
 \left( f (\theta)r^{-n+2k}\right)
^\wedge(x)&=&   (-1)^{k} \int_{S^{n-1}} \ln|(x,\xi)|\Delta^k
(f(\theta)r^{-n+2k})(\xi)d\xi\\
&+&  (-1)^{k} (2-n) \int_{S^{n-1}}  \Delta^{k-1} (f(\theta)r^{-n+2k})(\xi)d\xi.
\end{eqnarray*}
\end{enumerate}
\end{Lem}
\noindent{\bf Proof.} (i) and (ii) are proved in \cite[Lemma 3.16]{K-book}.

(iii)  is essentially from  \cite[Lemma 2.4]{Ya}. For completeness we include a proof.  For $q$
close to 1 we use part (i) with $k=1$ to get
\begin{equation}\label{eqn:limit}
\left( f (\theta)r^{-n+q+1}\right) ^\wedge(x)=\frac{\pi
  \int_{S^{n-1}} |(x,\xi)|^{1-q}\Delta
(f(\theta)r^{-n+q+1})(\xi)d\xi}{2
\Gamma(2-q)\sin\frac{\pi(1-q)}{2}}.
\end{equation}
When $q$ approaches $1$, both the numerator and denominator in the right hand side tend to zero.
Indeed, let us show that the limit of the numerator is zero:
$$\lim_{q\to 1}\int_{S^{n-1}} |(x,\xi)|^{1-q}\Delta
(f(\theta)r^{-n+q+1})(\xi)d\xi=\int_{S^{n-1}}\Delta (f(\theta)  r^{-n+2})(\xi)d\xi.$$ Recall the
relation between the spherical Laplacian $\Delta_S$ and Euclidean Laplacian $\Delta$ (see e.g.
\cite[p. 7]{Gr}). If $f$ is a homogeneous function of degree $m$, then on the sphere
$$\Delta_S f =\Delta f  -m(m+n-2)f .$$
Since $f(\theta)r^{-n+2}$ is a homogeneous function of degree $-n+2$, the previous formula implies
$\Delta (f(\theta)r^{-n+2})(\xi)= \Delta_S(f (\theta)r^{-n+2})(\xi)$. Due to the fact that
$\Delta_S$ is a self-adjoint operator, \cite[p. 7]{Gr}, we have
$$\int_{S^{n-1}}\Delta_S (f(\theta)r^{-n+2})(\xi)d\xi=0.$$

In order to compute the limit of (\ref{eqn:limit}) as $q \to 0$, apply l'Hopital's rule:
\begin{eqnarray*}
 \left( f (\theta)r^{-n+2}\right) ^\wedge(x)&=& -
 \int_{S^{n-1}} \ln|(x,\xi)|\Delta (f(\theta)r^{-n+2})(\xi)d\xi\\
&& -   \int_{S^{n-1}}  \Delta (f(\theta)r^{-n+2}\ln r)(\xi)d\xi.
\end{eqnarray*}
Computing the Laplacian in the latter integral and using Euler's formula for derivatives of
homogeneous functions, we get
\begin{eqnarray*}
 \left( f (\theta)r^{-n+2}\right) ^\wedge(x)&=& -
  \int_{S^{n-1}} \ln|(x,\xi)|\Delta (f(\theta)r^{-n+2})(\xi)d\xi\\
&& -   (2-n) \int_{S^{n-1}}   f(\xi) d\xi.
\end{eqnarray*}

Using the relation between the Fourier transform and differentiation, and applying the latter
formula to the function $\Delta^{k-1}(f (\theta)r^{-n+2k})$ which is  homogeneous  of degree
$-n+2$, we have
$$ \left( f (\theta)r^{-n+2k}\right) ^\wedge(x) = (-1)^{k-1} \left( \Delta^{k-1}(f
(\theta)r^{-n+2k})\right) ^\wedge(x)$$ $$ =   (-1)^{k} \int_{S^{n-1}} \ln|(x,\xi)|\Delta^k
(f(\theta)r^{-n+2k})(\xi)d\xi+$$
$$+  (-1)^{k} (2-n) \int_{S^{n-1}}  \Delta^{k-1}
(f(\theta)r^{-n+2k})(\xi)d\xi.$$

\qed

We will need the following spherical version of Parseval's  formula, for the proof see
\cite[Section 3.4]{K-book}.

\begin{Lem}
Let $K$ and $L$ be origin-symmetric infinitely smooth star bodies in $\mathbb{R}^n$ and $0<p<n$.
Then\index{spherical Parseval's formula}
\begin{equation}\label{eqn:Parseval}
\int_{S^{n-1}} (\|x\|_K^{-p} )^\wedge (\xi)  (\|x\|_L^{-n+p} )^\wedge (\xi)d\xi= (2\pi)^n
\int_{S^{n-1}} \|x\|_K^{-p}\|x\|_L^{-n+p} dx.
\end{equation}
\end{Lem}

In what follows $C$ will always be a non-zero constant, not necessarily the same in different
lines. We also use the notation $a(\epsilon)\sim b(\epsilon)$, meaning that $\displaystyle
\lim_{\epsilon\to 0} {a(\epsilon)}/{b(\epsilon)}=1$.

\begin{Lem}\label{Lem:order} Let $p$, $q>0$ be integers, $p+q\le n-2$.
\begin{enumerate}

\item[(i)] If $n-p-q-1$ is even, then
 for all $\xi \in S^{n-1}$
$$(|x|_2^{-q} \|x\|_E^{-p})^\wedge (\xi) \le C \epsilon^{-n+p+q+1}.$$

\item[(ii)] If $n-p-q-1$ is odd, then for every small $\alpha>0$ there exists a constant
$C_\alpha$, such that for all $\xi \in S^{n-1}$,
$$(|x|_2^{-q} \|x\|_E^{- p})^\wedge(\xi) \le C_\alpha \epsilon^{-n+p+q+1/(1+\alpha)}.$$

\item[(iii)] Moreover, in both cases,
$$(|x|_2^{-q} \|x\|_E^{-p})^\wedge(e_n)\sim C\epsilon^{-n+p+q+1}  .$$

\end{enumerate}

\end{Lem}

\noindent{\bf Proof.} (i) Let $n-p-q-1=2k$.  By Lemma \ref{Lem:K&M} we have
\begin{equation}\label{even}
(|x|_2^{-q} \|x\|_E^{-p})^\wedge (\xi)= \pi (-1)^k \int_{S^{n-1}\cap\xi^\bot} \Delta^k (|x|_2^{-q}
\|x\|_E^{-p}) dx.
\end{equation}

Applying $\Delta$ under the integral, each time we get a factor of $1/\epsilon^2$. This gives
$\epsilon^{-2k}=\epsilon^{-n+p+q+1}$. Finally use that $\|x\|_E^{-1}\le 1$ for $x\in S^{n-1}$.


To prove (iii) for this case, we use (\ref{even}) again. Keeping only the terms of the highest
order of $\epsilon$, we get
$$(|x|_2^{-q} \|x\|_E^{-p})^\wedge (e_n)\sim C \int_{S^{n-1}\cap
e_n^\bot}\frac{\partial^{2k}}{\partial x_n^{2k}}   \|x\|_E^{-p} dx \sim C
\frac{\partial^{2k}}{\partial x_n^{2k}} (1+\frac{x_n^2}{\epsilon^2})^{-p}\Big|_{x_n=0} = C
\epsilon^{-2k}. $$

(ii) Let  $n-p-q-1=2k-1$. By Lemma \ref{Lem:K&M},
\begin{eqnarray}\label{odd}
(|x|_2^{-q} \|x\|_E^{-p})^\wedge (\xi)&=& (-1)^k \int_{S^{n-1}} \ln|(x,\xi)| \Delta^k (|x|_2^{-q}
\|x\|_E^{-p}) dx \nonumber \\
&+&(-1)^k \int_{S^{n-1}} \Delta^{k-1} (|x|_2^{-q} \|x\|_E^{-p}) dx.
\end{eqnarray}

Consider the first integral in (\ref{odd}).  As before, applying the Laplacian $k$ times under the
integral, we get a factor of $\epsilon^{-2k}=\epsilon^{-n+p+q}$. Therefore we need to estimate
terms of the following form
$$ \int_{S^{n-1}}   \ln
|(\theta,\xi)  \|\theta\|_E^{-n+q} d\theta.$$

By H\"older's inequality, for a small $\alpha >0$,
$$\left|\int_{S^{n-1}}   \ln
|(\theta,\xi)  \|\theta\|_E^{-n+q} d\theta \right|$$
$$\le  \left(\int_{S^{n-1}}  \Big|\ln |(\theta,\xi)|\Big|^{(1+\alpha)/\alpha} d\theta
\right)^{\alpha/(1+\alpha)} \left( \int_{S^{n-1}} \|\theta\|_E^{(-n+q)(1+\alpha)}
d\theta\right)^{1/(1+\alpha)}$$
$$= C_\alpha \left( \int_{S^{n-1}} \|\theta\|_E^{(-n+q)(1+\alpha)}
d\theta\right)^{1/(1+\alpha)}.$$

Note that Parseval's formula (\ref{eqn:Parseval})  gives
\begin{eqnarray}\label{integral}
&&\int_{S^{n-1}} \|\theta\|_E^{(-n+q)(1+\alpha)} d\theta \nonumber\\
&&= (2\pi)^{-n}\int_{S^{n-1}} \left(\|x\|_E^{(-n+q)(1+\alpha)}\right)^\wedge(\xi) \left(|x|_2^{
n\alpha -q-q\alpha}\right)^\wedge(\xi) d\xi \nonumber\\
&& =C \epsilon \int_{S^{n-1}} \|\theta\|_{E^*}^{n\alpha -q-q\alpha} d\theta,
\end{eqnarray}
where $E^*$ is the ellipsoid given by
\begin{equation}\label{E^*}\|x\|_{E^*}=  \left( x_1^2+\cdots
+ x_{n-1}^2 +{\epsilon^2}{ x_n^2} \right)^{1/2}.
\end{equation}

By the following elementary formula (see e.g. \cite[p.9]{Gr})
\begin{equation}\label{formula:t-coord}
\int_{S^{n-1}}f((x,\theta))d\theta= |S^{n-2}| \int_{-1}^{1} (1-t^2)^{(n-3)/2} f(t)dt, \quad   x\in
S^{n-1},
\end{equation}
the integral in (\ref{integral}) equals
$$=C \epsilon
\int_{-1}^1 (1-t^2)^{(n-3)/2}  (1-t^2+\epsilon^2t^2)^{(n\alpha -q-q\alpha)/2} dt$$ $$\sim C
\epsilon \int_{-1}^1 (1-t^2)^{(n-3+n\alpha -q-q\alpha)/2} dt = C \epsilon. $$ The latter integral
converges if $(n-q)(1+\alpha)>1$, which is the case.

Therefore for all small $\alpha>0$,
$$\left|\int_{S^{n-1}}   \ln
|(\theta,\xi)  \|\theta\|_E^{-n+q} d\theta \right|\le C_\alpha \epsilon^{1/(1+\alpha)},$$ and hence
the first integral in (\ref{odd}) can be bounded as follows:
$$\left|\int_{S^{n-1}} \ln|(x,\xi)| \Delta^k (|x|_2^{-q} \|x\|_E^{-p}) dx\right|\le
C_\alpha \epsilon^{-n+p+q+1/(1+\alpha)}.$$

We use similar ideas to estimate the second integral in (\ref{odd}). Applying Parseval's formula
two times and using the relation between the Fourier transform and differentiation, we get
$$\int_{S^{n-1}} \Delta^{k-1} (|x|_2^{-q} \|x\|_E^{-p}) dx= \int_{S^{n-1}} \Delta^\frac{n-p-q-2 }{2} (|x|_2^{-q} \|x\|_E^{-p}) |x|_2^{-2} dx$$
$$ = (2\pi)^{-n} \int_{S^{n-1}} \left(\Delta^\frac{n-p-q-2 }{2} (|x|_2^{-q} \|x\|_E^{-p})\right)^\wedge(\theta) \left(|x|_2^{-2}\right)^\wedge(\theta) d\theta$$
$$ = C \int_{S^{n-1}} \left( |x|_2^{-q} \|x\|_E^{-p}\right)^\wedge(\theta)   d\theta  = C \int_{S^{n-1}}   \|x\|_E^{-p}   d\theta =O(1),$$
since $\|x\|_E^{-1}\le 1$ on the sphere.

Combining the estimates for the both terms in (\ref{odd}), one can see that
$$(|x|_2^{-q} \|x\|_E^{- q})^\wedge(\xi) \le C_\alpha \epsilon^{-n+p+q+1/(1+\alpha)}.$$

Now we will show that almost the same degree of dependence is achieved when $\xi=e_n$.   Indeed,
using formula (\ref{odd}) with $\xi=e_n$ and dropping the second integral (which is small compared
to the first integral), we get
$$(|x|_2^{-p} \|x\|_E^{-q})^\wedge(e_n)\sim (-1)^k \int_{S^{n-1}}
\ln|x_n| \Delta^k (|x|_2^{-q} \|x\|_E^{-p}) dx.
$$

Formula (\ref{formula:t-coord}) applied to the last integral gives
$$=2(-1)^k \int_{0}^1 \ln x_n \cdot (1-x_n^2)^{\frac{n-3}{2}}
  \Delta^k (|x|_2^{-q} \|x\|_E^{-p})\Big|_{x_1^2+\cdots+x_{n-1}^2=1-x_n^2} dx_n.
$$
After the change of the variable $x_n=\epsilon\cdot z$, the latter equals
$$=2(-1)^k \epsilon \int_{0}^{1/\epsilon} \ln (\epsilon z) \cdot (1-\epsilon^2 z^2)^{\frac{n-3}{2}}
   \Delta_\epsilon^k (\|x\|_{E^*}^{-q}
|x|_2^{-p})\Big|_{ x_1^2+\cdots+x_{n-1}^2=1- \epsilon^2 z^2, x_n=  z }
  dz,
$$
where $E^*$ is defined by (\ref{E^*}), and
 $$\Delta_\epsilon = \frac{\partial^2}{\partial x_1^2}+\cdots + \frac{\partial^2}{\partial
x_{n-1}^2}+\frac{1}{\epsilon^2 } \frac{\partial^2}{  \partial x_n^2}.$$

Note that after writing $\ln(\epsilon z)=\ln \epsilon + \ln z$, we will have two integrals, the
first being equal to
$$2(-1)^k \epsilon \ln \epsilon  \int_{0}^{1/\epsilon}   (1-\epsilon^2
z^2)^{\frac{n-3}{2}}
  \Delta_\epsilon^k (\|x\|_{E^*}^{-q} |x|_2^{-p})\Big|_{x_1^2+\cdots+x_{n-1}^2=1- \epsilon^2 z^2, x_n=  z}
  dz
$$
$$=(-1)^k \ln \epsilon  \int_{S^{n-1}}
  \Delta^k (|x|_2^{-q} \|x\|_E^{-p}) dx =(-1)^k \ln \epsilon  \int_{S^{n-1}}
 \Delta (\Delta^{k-1} (|x|_2^{-q} \|x\|_E^{-p})) dx.
$$

Under the integral we have the Laplacian of a homogeneous function of degree $-n+2$ which equals
the spherical Laplacian. Since the spherical Laplacian is a self-adjoint operator,  the latter
integral is equal to zero.

Therefore we only need to compute the order of the integral
$$2(-1)^k \epsilon \int_{0}^{1/\epsilon} \ln  z \cdot (1-\epsilon^2 z^2)^{\frac{n-3}{2}}
  \Delta_\epsilon^k (\|x\|_{E^*}^{-q} |x|^{-p})\Big|_{x_1^2+\cdots+x_{n-1}^2=1- \epsilon^2 z^2,x_n=  z}
  dz.
$$

As before, the largest term is obtained when we apply $\frac{1}{\epsilon^2 } \frac{\partial^2}{
\partial x_n^2}$
to $|x|^{-p}$ successively  $k$ times.
$$\sim 2(-1)^k   {\epsilon^{-2k+1} }\int_{0}^{1/\epsilon} \ln  z \cdot (1-\epsilon^2
z^2)^{\frac{n-3}{2}}
    \frac{\partial^{2k }}{
\partial x_n^{2k }} (|x|^{-p})\Big|_{x_1^2+\cdots+x_{n-1}^2=1- \epsilon^2 z^2,x_n=  z}
  dz.
$$

It is enough to show that $$\int_{0}^{1/\epsilon} \ln  z \cdot (1-\epsilon^2 z^2)^{\frac{n-3}{2}}
    \frac{\partial^{2k}}{
\partial x_n^{2k}} (|x|^{-p})\Big|_{x_1^2+\cdots+x_{n-1}^2=1-  \epsilon^2z^2,x_n=   z}
  dz
$$ has a finite nonzero limit as $\epsilon\to 0$.

One can see that
$$\lim_{\epsilon\to 0} \int_{0}^{1/\epsilon} \ln  z \cdot (1-\epsilon^2 z^2)^{\frac{n-3}{2}}
    \frac{\partial^{2k }}{
\partial x_n^{2k }} (|x|^{-p})\Big|_{x_1^2+\cdots+x_{n-1}^2=1-  \epsilon^2z^2,x_n=   z}
  dz
$$
$$= \int_{0}^{\infty} \ln  z \cdot
    \frac{\partial^{2k }}{
\partial z^{2k}} (1+z^2)^{-p/2}   dz.
$$

To finish the proof, we need to show that the latter integral is not equal to zero. Let
$P(z^2)=c_0+c_1 z^2+\cdots+c_{k-1} z^{2k-2}$ be the Taylor polynomial of $(1+z^2)^{-p/2}$ at zero
of order $2k-2$. Then clearly,
$$\int_{0}^{\infty} \ln  z \    \frac{\partial^{2k }}{ \partial z^{2k }} (1+z^2)^{-p/2}   dz
 =\int_{0}^{\infty} \ln  z \     \frac{\partial^{2k }}{ \partial z^{2k }}\left( (1+z^2)^{-p/2}-P(z^2)\right)
 dz.
$$
After integration by parts $2k$ times and the change of the variable $t=z^2$ the integral becomes
$$= - (2k-1)!\int_{0}^{\infty}   z^{-2k}  \left( (1+z^2)^{-p/2}-P(z^2)\right)   dz
$$
$$= -\frac12 (2k-1)!\int_{0}^{\infty}   t^{-k-1/2}  \left( (1+t)^{-p/2}-P(t)\right)   dt.
$$
Using integration by parts in the opposite order and observing that $P(t)$ is the Taylor polynomial
of $(1+t)^{-p/2}$, we get
$$= -\frac12 \frac{(2k-1)!}{(1/2)(3/2)\cdots (k-1/2)}\int_{0}^{\infty}   t^{-1/2} \frac{\partial^{k }}{ \partial t^{k }} \left( (1+t)^{-p/2}-P(t)\right)   dt
$$
$$= -\frac12 \frac{(2k-1)!}{(1/2)(3/2)\cdots (k-1/2)}\int_{0}^{\infty}   t^{-1/2} \frac{\partial^{k }}{ \partial t^{k }} \left( (1+t)^{-p/2}\right)
dt.
$$

The latter is clearly a nonzero constant.

 \qed

 Now we are ready to prove the main result of this section.

\begin{Thm}
For every $1\le k< l < n-3$ there exists an origin-symmetric convex body $K\in \mathbb{R}^n$ that
does not belong to $\mathcal I_k$, but belongs to $\mathcal I_{l}$.
\end{Thm}
\noindent{\bf Proof.}  For a small $\epsilon>0$   define a body $K$ by
\begin{eqnarray*}\label{Def:K}
\|x\|_K^{-1}=|x|_2^{-1} -     \epsilon^{n-k-3/2} \|x\|_E^{-1}, \quad x \in \mathbb{R}^n
\setminus\{0\},
\end{eqnarray*}
where   $E$ is the  ellipsoid with the norm
$$\|x\|_E^{-1}=  \left( x_1^2+\cdots +
x_{n-1}^2 +\frac{ x_n^2}{\epsilon^2}\right)^{-1/2}.$$

Convexity of $K$ follows along the lines of Lemma \ref{K-convex}. (One will need
$\epsilon^{n-k-3/2-2}$ to be small, which is the case since $n-3>k$).

Consider the $-l$th power of the norm of $K$.
$$\|x\|_K^{-l}=|x|_2^{-l} -    \epsilon^{n-k-3/2} l |x|_2^{-l+1} \|x\|_E^{-1} + R(x),$$
where $R(x)$ is the sum of the terms of the following form:
$$\epsilon^{i(n-k-3/2)}  |x|_2^{-l+i} \|x\|_E^{-i}, \qquad i\ge 2.$$

Applying the Fourier transform, we have for all $\xi\in S^{n-1}$,
$$(\|x\|_K^{-l})^\wedge (\xi) =C(l,n) -    \epsilon^{n-k-3/2} l (|x|_2^{-l+1} \|x\|_E^{-1})^\wedge(\xi) + \widehat{R}(\xi).$$
By Lemma \ref{Lem:order} the order of the second term is at most
$$\epsilon^{n-k-3/2} \epsilon^{-n+l+1/(1+\alpha)}=\epsilon^{l-k-1/2-\alpha/(1+\alpha)}\to 0, \mbox{ as } \epsilon \to 0,$$
if $\alpha$ is small enough. $\widehat{R}$ is even smaller, since it contains terms of the order
$\epsilon^{i(n-k-3/2)} \epsilon^{-n+l+1/(1+\alpha)},$ $ i\ge 2.$

Therefore if $\epsilon$ is small, then $(\|x\|_K^{-l})^\wedge (\xi)\ge 0$, and so $K\in \mathcal
I_{l}$.

Now consider the $-k$th power of the norm of $K$.
$$\|x\|_K^{-k}=|x|_2^{-k} -    \epsilon^{n-k-3/2} k |x|_2^{-k+1} \|x\|_E^{-1} + Q(x).$$

Computing the Fourier transform in the direction of $\xi=e_n$, we have
\begin{equation}\label{-k}
(\|x\|_K^{-k})^\wedge (e_n) =C(k,n) -    \epsilon^{n-k-3/2} k (|x|_2^{-k+1}
\|x\|_E^{-1})^\wedge(e_n) + \widehat{Q}(e_n).
\end{equation}
$\widehat{Q}(e_n)$ is small, since it has terms of order at most $\epsilon^{i(n-k-3/2)}
\epsilon^{-n+k+1}$, $ i\ge 2.$ Therefore, we will pay attention only to the second term in
(\ref{-k}). By Lemma \ref{Lem:order},
$$ \epsilon^{n-k-3/2}(|x|_2^{-k+1} \|x\|_E^{-1})^\wedge(e_n) \sim C \epsilon^{n-k-3/2}\epsilon^{-n+k+1} = C\epsilon^{-1/2}.$$
If we choose $\epsilon>0$ small enough so that the latter is greater than $C(k,n)$, then
$(\|x\|_K^{-k})^\wedge (e_n)<0$. So $K\not\in \mathcal I_k$.

\qed

{\bf Acknowledgements.} The author is thankful to Professors Alexander Koldobsky and Paul Goodey
for many valuable suggestions.

\end{document}